\documentclass[12pt]{amsart}
\usepackage[utf8]{inputenc}
\usepackage{amsmath,amssymb,amsthm}
\usepackage{enumerate, euler, palatino}
\usepackage[colorlinks=true]{hyperref}
\hypersetup{urlcolor=blue, citecolor=red , linkcolor= blue}
\usepackage{cite}
\usepackage[a4paper,margin=2.5cm]{geometry}

\theoremstyle{plain}
\newtheorem{theorem}{Theorem}[section]
\newtheorem{proposition}[theorem]{Proposition}

\newtheorem{corollary}[theorem]{Corollary}
\theoremstyle{definition}

\newtheorem{remark}[theorem]{Remark}
\newtheorem{question}[theorem]{Question}

\title[Quasinilpotency of commutators]{Quasinilpotency of commutators under a weakened Kleinecke--Shirokov condition}
\author{Hranislav Stankovi\'c}
\address{Faculty of Electronic Engineering, University of Ni\v{s}, Aleksandra Medvedeva 4, 18000 Ni\v{s}, Serbia}
\email{hranislav.stankovic@elfak.ni.ac.rs}

\begin{document}
	
	\begin{abstract}
		Kramar and Kramar Fijav\v{z} proved that a commutator $Y=AB-BA$ of two complex square matrices is nilpotent whenever $Y$ commutes with an iterated commutator $\delta_Y^{k}(A)$, where $\delta_Y(T)=YT-TY$, and asked whether this result can be generalized to bounded operators on infinite-dimensional spaces. We show that for $k\geq 2$ such a commutator need not be quasinilpotent. On the other hand, for every $k$ it is quasinilpotent whenever it is compact, or more generally whenever it is a Riesz operator or $AY-YA$ is compact. Moreover, its Riesz projection at any non-zero isolated point of its spectrum always has infinite rank.
	\end{abstract}
	
	\maketitle
	\bigskip
	\section{Introduction}
	\bigskip
	
	Let $\mathcal{A}$ be a complex unital Banach algebra. For $a\in\mathcal{A}$ we denote by $\delta_a$ the inner derivation
	\[
	\delta_a(t)=at-ta,\qquad t\in\mathcal{A},
	\]
	and we write $[a,b]=\delta_a(b)=ab-ba$. The Kleinecke--Shirokov theorem \cite{Kleinecke,Shirokov} asserts that if $a,b\in\mathcal{A}$ are such that $a$ commutes with $[a,b]$, then $[a,b]$ is quasinilpotent. Its finite-dimensional predecessor, Jacobson's lemma \cite{Jacobson}, says that if $A$ and $B$ are complex square matrices and $A$ commutes with $AB-BA$, then $AB-BA$ is nilpotent. For localizations and other variations of the Kleinecke--Shirokov theorem we refer to \cite{BracicKuzma,Kramar2016,Stankovic2025,Stankovic2026a}.
	
	Recently, Kramar and Kramar Fijav\v{z} \cite{KramarKramar} weakened the commutativity hypothesis in Jacobson's lemma. Writing $Y=AB-BA$, they proved the following two results for complex square matrices $A$ and $B$.
	\begin{enumerate}[\rm($i$)]
		\item If $Y$ commutes with $\delta_A(Y)=AY-YA$, then $Y$ is nilpotent \cite[Theorem 1]{KramarKramar}.
		\item If $Y$ commutes with $\delta_Y^{k}(A)$ for some $k\in\mathbb{N}$, then $Y$ is nilpotent \cite[Theorem 2]{KramarKramar}.
	\end{enumerate}
	Since $\delta_A(Y)=-\delta_Y(A)$, the hypothesis in ($i$) is the case $k=1$ of that in ($ii$), and both can be written as
	\begin{equation}\label{eq:hyp}
		\delta_Y^{k+1}(A)=0 .
	\end{equation}
	In the same paper the authors posed the following problem \cite[Question 1]{KramarKramar}: can these two results be generalized to bounded operators on infinite-dimensional spaces? In that setting nilpotency is naturally replaced by quasinilpotency, exactly as in the Klei\-necke--Shirokov theorem, so the problem reads as follows.
	
	\begin{question}\label{q:main}
		Let $\mathcal{X}$ be an infinite-dimensional complex Banach space, let $A,B\in\mathfrak{B}(\mathcal{X})$, $Y=AB-BA$, and suppose that $\delta_Y^{k+1}(A)=0$ for some $k\in\mathbb{N}$. Is $Y$ quasinilpotent?
	\end{question}
	
	In Section \ref{sec:negative} we show that for $k\geq 2$ the answer is negative. Using the Brown--Pearcy characterization of commutators in $\mathfrak{B}(\mathcal{H})$, we exhibit, on any separable infinite-dimensional Hilbert space $\mathcal{H}$, a commutator $Y=AB-BA$ with $\sigma(Y)=\{1\}$ such that $\delta_Y^{3}(T)=0$ for every $T\in\mathfrak{B}(\mathcal{H})$.
	
	In Section \ref{sec:riesz} we use the fact that, under the hypothesis \eqref{eq:hyp}, the element $a$ commutes with every Riesz idempotent $p$ of $y=[a,b]$, so that the hypothesis passes to the commutator $yp=[pap,pbp]$ in the algebra $p\mathcal{A}p$. For operators this yields, by a trace argument, that every non-zero isolated point of $\sigma(Y)$ has an infinite-dimensional Riesz subspace. In particular, the answer to Question \ref{q:main} is affirmative, for every $k$, whenever $Y$ is compact, or more generally whenever $Y$ is a Riesz operator or $\delta_A(Y)$ is compact.
	
	\bigskip 
	
	Let us fix the notation and recall the facts we shall use. For a complex Banach space $\mathcal{X}$ we write $\mathfrak{B}(\mathcal{X})$ for the algebra of all bounded linear operators on $\mathcal{X}$ and $\mathcal{K}(\mathcal{X})$ for its ideal of compact operators. The spectrum and the spectral radius of an element $y$ of a Banach algebra are denoted by $\sigma(y)$ and $r(y)$, and $y$ is quasinilpotent if $r(y)=0$. Since $\delta_a$ is the difference of the operators of left and right multiplication by $a$, and these two operators commute, the binomial formula gives
	\begin{equation}\label{eq:binom}
		\delta_a^{n}(t)=\sum_{j=0}^{n}(-1)^j\binom{n}{j}a^{n-j}ta^{j},\qquad t\in\mathcal{A},\ n\in\mathbb{N}.
	\end{equation}
	
	Let $y\in\mathcal{A}$ and let $\sigma_1$ be a spectral set of $y$, that is, a subset of $\sigma(y)$ which is both open and closed in $\sigma(y)$. The Riesz idempotent of $y$ associated with $\sigma_1$ is
	\[
	p=\frac{1}{2\pi i}\int_{\Gamma}(z-y)^{-1}\,dz ,
	\]
	where $\Gamma$ is a contour surrounding $\sigma_1$ and separating it from $\sigma(y)\setminus\sigma_1$. We recall that $p$ is an idempotent commuting with $y$, that $yp$ belongs to the unital Banach algebra $p\mathcal{A}p$ with unit $p$ and has spectrum $\sigma_1$ there, and that $y(1-p)$ has spectrum $\sigma(y)\setminus\sigma_1$ in $(1-p)\mathcal{A}(1-p)$. If $y=Y\in\mathfrak{B}(\mathcal{X})$ and $\sigma_1=\{\lambda\}$, where $\lambda$ is an isolated point of $\sigma(Y)$, the range of the corresponding Riesz projection is called the Riesz subspace of $Y$ at $\lambda$, and $\lambda$ is called a Riesz point of $Y$ if this subspace is finite-dimensional. An operator $Y\in\mathfrak{B}(\mathcal{X})$ is a Riesz operator if every non-zero point of $\sigma(Y)$ is a Riesz point, and by Ruston's characterization \cite{Ruston} this holds if and only if the image of $Y$ in the Calkin algebra $\mathfrak{B}(\mathcal{X})/\mathcal{K}(\mathcal{X})$ is quasinilpotent (see also \cite{Aiena}).
	
	We shall use the following classical result of Lumer and Rosenblum \cite{LumerRosenblum} (see also \cite{Rosenblum}): if $L$ and $R$ are commuting bounded operators on a Banach space, then $\sigma(L-R)\subseteq\sigma(L)-\sigma(R)$. Finally, we recall the theorem of Brown and Pearcy \cite{BrownPearcy}: if $\mathcal{H}$ is a separable infinite-dimensional Hilbert space, then an operator $T\in\mathfrak{B}(\mathcal{H})$ is a commutator, that is, $T=AB-BA$ for some $A,B\in\mathfrak{B}(\mathcal{H})$, if and only if $T$ is not of the form $\lambda I+K$ with $\lambda\neq 0$ and $K\in\mathcal{K}(\mathcal{H})$.
	
	\bigskip
	\section{A counterexample for \texorpdfstring{$k\geq 2$}{k >= 2}}\label{sec:negative}
	\bigskip
	
	The following proposition shows that the answer to Question \ref{q:main} is negative for every $k\geq 2$.
	
	\begin{proposition}\label{prop:negative}
		Let $\mathcal{H}$ be a separable infinite-dimensional complex Hilbert space. There exist $A,B\in\mathfrak{B}(\mathcal{H})$ such that the commutator $Y=AB-BA$ satisfies
		\[
		\delta_Y^{3}(T)=0\quad\text{for every }T\in\mathfrak{B}(\mathcal{H})
		\]
		and $\sigma(Y)=\{1\}$. In particular, $\delta_Y^{k+1}(A)=0$ for every $k\geq 2$, while $Y$ is not quasinilpotent.
	\end{proposition}
	
	\begin{proof} 
		Write $\mathcal{H}=\mathcal{H}_0\oplus\mathcal{H}_0$, where $\mathcal{H}_0$ is a separable infinite-dimensional Hilbert space, and consider
		\[
		N=\begin{pmatrix} 0 & I\\ 0 & 0\end{pmatrix},\qquad Y=I+N=\begin{pmatrix} I & I\\ 0 & I\end{pmatrix}.
		\]
		Then $N^2=0$, so that $\sigma(Y)=\{1\}$. We claim that $Y$ is a commutator. By the Brown--Pearcy theorem it suffices to check that $Y$ is not of the form $\lambda I+K$ with $\lambda\neq 0$ and $K$ compact. Suppose that $I+N=\lambda I+K$ with $K$ compact. Then $N=(\lambda-1)I+K$, and squaring gives
		\[
		0=N^2=(\lambda-1)^2I+2(\lambda-1)K+K^2 ,
		\]
		so that $(\lambda-1)^2I=-2(\lambda-1)K-K^2$ is compact. Since $\mathcal{H}$ is infinite-dimensional, $\lambda=1$, and hence $N=K$ is compact. This is impossible, since the compactness of $N$ would imply the compactness of the identity operator on the infinite-dimensional space $\mathcal{H}_0$. Consequently there exist $A,B\in\mathfrak{B}(\mathcal{H})$ with $AB-BA=Y$.
		
		Since $\delta_Y=\delta_{I+N}=\delta_N$, formula \eqref{eq:binom} gives, for every $T\in\mathfrak{B}(\mathcal{H})$,
		\[
		\delta_Y^{3}(T)=\delta_N^{3}(T)=N^3T-3N^2TN+3NTN^2-TN^3=0 ,
		\]
		because $N^2=0$.
	\end{proof}
	
	\begin{remark}
		The construction above does not settle the case $k=1$. Indeed, for $Y=I+N$ as in the proof, formula \eqref{eq:binom} gives
		\[
		\delta_Y^{2}(A)=\delta_N^{2}(A)=N^2A-2NAN+AN^2=-2NAN ,
		\]
		so that $\delta_Y^{2}(A)=0$ if and only if $NAN=0$. Thus $Y=I+N$ gives a negative answer to Question \ref{q:main} for $k=1$ if and only if there exist $A,B\in\mathfrak{B}(\mathcal{H})$ such that $AB-BA=I+N$ and $NAN=0$.
	\end{remark}
	
	\bigskip
	\section{Riesz points of the commutator}\label{sec:riesz}
	\bigskip
	
	We now turn to positive results, which are valid for every $k$ and rest on the Riesz decomposition of the commutator. For operators, assertion ($i$) of the following proposition is contained in \cite[Lemma 2.25(iii)]{LomonosovShulman} (see also \cite[Proposition 2.2]{LaursenNeumann}), and the general case reduces to this one by means of the left regular representation of $\mathcal{A}$. For the convenience of the reader we include a short direct proof.
	
	\begin{proposition}\label{prop:riesz}
		Let $a,y\in\mathcal{A}$ and suppose that $\delta_y^{k+1}(a)=0$ for some $k\in\mathbb{N}$. Let $\sigma_1$ be a spectral set of $y$ and let $p$ be the associated Riesz idempotent. Then the following hold.
		\begin{enumerate}[\rm($i$)]
			\item $ap=pa$.
			\item If $y=[a,b]$ for some $b\in\mathcal{A}$, then in the unital Banach algebra $p\mathcal{A}p$ one has $yp=[pap,pbp]$ and moreover $\delta_{yp}^{k+1}(pap)=0$.
		\end{enumerate}
	\end{proposition}
	
	\begin{proof}
		Put $q=1-p$. Since $p$ and $q$ commute with $y$, for every $t\in\mathcal{A}$ we have
		\[
		\delta_y(ptq)=yptq-ptqy=pytq-ptyq=p\,\delta_y(t)\,q .
		\]
		In particular, the subspace $p\mathcal{A}q$ is invariant under $\delta_y$. It is closed in $\mathcal{A}$, being the range of the bounded idempotent $t\mapsto ptq$, and hence a Banach space. Since $p$ and $q$ are idempotents, we also have $yptq=(yp)(ptq)$ and $ptqy=(ptq)(yq)$, so that the restriction of $\delta_y$ to $p\mathcal{A}q$ equals $L-R$, where $L$ and $R$ denote the operators of left multiplication by $yp$ and of right multiplication by $yq$ on the Banach space $p\mathcal{A}q$. The operators $L$ and $R$ commute. If $\lambda\notin\sigma_1$, then $\lambda p-yp$ is invertible in $p\mathcal{A}p$, and we denote its inverse by $u\in p\mathcal{A}p$. Let $U$ be the operator of left multiplication by $u$ on $p\mathcal{A}q$, that is, $U(ptq)=uptq$. Since $(\lambda-L)(ptq)=(\lambda p-yp)ptq$ for every $t\in\mathcal{A}$, we obtain
		\[
		U(\lambda-L)(ptq)=ptq\quad\text{and}\quad(\lambda-L)U(ptq)=ptq .
		\]
		Thus $\lambda-L$ is invertible on $p\mathcal{A}q$ with inverse $U$, and therefore $\sigma(L)\subseteq\sigma_1$. In the same way, if $\lambda\notin\sigma(y)\setminus\sigma_1$, then the operator of right multiplication by the inverse of $\lambda q-yq$ in $q\mathcal{A}q$ is the inverse of $\lambda-R$ on $p\mathcal{A}q$, so that $\sigma(R)\subseteq\sigma(y)\setminus\sigma_1$. By the Lumer--Rosenblum theorem,
		\[
		\sigma(L-R)\subseteq\sigma_1-\big(\sigma(y)\setminus\sigma_1\big),
		\]
		and the set on the right does not contain $0$. Therefore the restriction of $\delta_y$ to $p\mathcal{A}q$ is invertible.
		
		Iterating the identity $\delta_y(ptq)=p\,\delta_y(t)\,q$, we get $\delta_y^{k+1}(paq)=p\,\delta_y^{k+1}(a)\,q=0$, and the invertibility of $\delta_y$ on $p\mathcal{A}q$ gives $paq=0$. Replacing $p$ by $q$ in the argument above we obtain $qap=0$ as well. Hence
		\[
		ap=(p+q)ap=pap=pa(p+q)=pa ,
		\]
		which proves ($i$).
		
		For ($ii$), let $y=[a,b]$. Since $p$ commutes with $y$ and $p^2=p$, we have $yp=yp^2=pyp$. Using also $ap=pa$ and $p^2=p$, we obtain
		\begin{align*}
			yp&=pyp=p(ab-ba)p=pabp-pbap\\
			&=papbp-pbpap=(pap)(pbp)-(pbp)(pap)=[pap,pbp].
		\end{align*}
		Finally, since $p$ commutes with $y$ and $p^2=p$, for every $t\in\mathcal{A}$ we have
		\[
		\delta_{yp}(ptp)=(yp)(ptp)-(ptp)(yp)=pytp-ptyp=p\,\delta_y(t)\,p .
		\]
		By induction on $n$ this gives $\delta_{yp}^{n}(ptp)=p\,\delta_y^{n}(t)\,p$ for all $n\in\mathbb{N}$ and $t\in\mathcal{A}$. Indeed, if this holds for some $n$, then
		\[
		\delta_{yp}^{n+1}(ptp)=\delta_{yp}\big(p\,\delta_y^{n}(t)\,p\big)=p\,\delta_y\big(\delta_y^{n}(t)\big)\,p=p\,\delta_y^{n+1}(t)\,p .
		\]
		Taking $t=a$ and $n=k+1$, we obtain
		\[
		\delta_{yp}^{k+1}(pap)=p\,\delta_y^{k+1}(a)\,p=0 .\qedhere
		\]
	\end{proof}
	
	We now draw the consequences for operators.
	
	\begin{corollary}\label{cor:isolated}
		Let $\mathcal{X}$ be a complex Banach space, $A,B\in\mathfrak{B}(\mathcal{X})$, $Y=AB-BA$, and suppose that $\delta_Y^{k+1}(A)=0$ for some $k\in\mathbb{N}$. Then the following hold.
		\begin{enumerate}[\rm($i$)]
			\item If $\lambda\in\sigma(Y)\setminus\{0\}$ is an isolated point of $\sigma(Y)$, then the Riesz subspace of $Y$ at $\lambda$ is infinite-dimensional.
			\item If $Y$ is a Riesz operator, in particular if $Y$ is compact, then $Y$ is quasinilpotent.
		\end{enumerate}
	\end{corollary}
	
	\begin{proof}
		Let $\lambda\neq 0$ be an isolated point of $\sigma(Y)$, let $P$ be the Riesz projection of $Y$ at $\lambda$ and put $M=P(\mathcal{X})$. By Proposition \ref{prop:riesz}, applied in the Banach algebra $\mathfrak{B}(\mathcal{X})$ with $p=P$, we have $YP=[PAP,PBP]$. For $x\in M$ we have $Px=x$, and hence $Yx=YPx=PYx\in M$, because $P$ commutes with $Y$. Clearly, $PAP$ and $PBP$ map $\mathcal{X}$ into $M$. Restricting the identity $YP=[PAP,PBP]$ to $M$, where $YP$ coincides with $Y$, we obtain $Y|_M=[(PAP)|_M,(PBP)|_M]$. If $M$ were finite-dimensional, taking traces would give $\operatorname{tr}(Y|_M)=0$. On the other hand, $\sigma(Y|_M)=\{\lambda\}$, so that $\operatorname{tr}(Y|_M)=\lambda\dim M\neq 0$. This contradiction proves ($i$). If $Y$ is a Riesz operator, every non-zero point of $\sigma(Y)$ is an isolated point with finite-dimensional Riesz subspace, so it follows from ($i$) that $\sigma(Y)=\{0\}$. This applies in particular to compact operators, which are Riesz operators by the Riesz--Schauder theory.
	\end{proof}
	
	Corollary \ref{cor:isolated} shows that both theorems of Kramar and Kramar Fijav\v{z} remain valid for compact commutators on infinite-dimensional Banach spaces, with nilpotency replaced by quasinilpotency. The compactness of $Y$ can be weakened to the compactness of $\delta_A(Y)=AY-YA$.
	
	\begin{corollary}\label{cor:derivation}
		Let $\mathcal{X}$ be a complex Banach space, $A,B\in\mathfrak{B}(\mathcal{X})$, $Y=AB-BA$, and suppose that $\delta_Y^{k+1}(A)=0$ for some $k\in\mathbb{N}$. If $\delta_A(Y)$ is compact, then $Y$ is quasinilpotent.
	\end{corollary}
	
	\begin{proof}
		Since $\delta_A(Y)$ is compact, $Y$ is a Riesz operator by a theorem of Lin \cite{Lin} (see also \cite[Theorem 3]{Kramar2016}). Indeed, if $\pi\colon\mathfrak{B}(\mathcal{X})\to\mathfrak{B}(\mathcal{X})/\mathcal{K}(\mathcal{X})$ is the canonical map, then $\pi(A)$ commutes with $\pi(Y)=[\pi(A),\pi(B)]$, so that $\pi(Y)$ is quasinilpotent by the Kleinecke--Shirokov theorem, and Ruston's characterization applies. Now Corollary \ref{cor:isolated} shows that $Y$ is quasinilpotent.
	\end{proof}
	
	\bigskip
	\section*{Declarations}
	\noindent{\bf{Funding}}\\
	This work has been supported by the Ministry of Science, Technological Development and Innovation of the Republic of Serbia [Grant Number: 451-03-34/2026-03/200102].
	
	\bigskip
	
	\noindent{\bf{Availability of data and materials}}\\
	\noindent No data were used to support this study.
	\bigskip\\
	\noindent{\bf{Competing interests}}\\
	\noindent The author declares that he has no competing interests.
	\bigskip\\
	\noindent{\bf{Declaration of generative AI and AI-assisted technologies in the manuscript preparation process}}\\
	\noindent During the preparation of this work, the author used ChatGPT (OpenAI) and Claude (Anthropic) to identify relevant references, to verify the mathematical arguments, to assist in the derivation and drafting of parts of the results, and for language editing and improvement of the exposition. The author reviewed, verified and edited all generated content and takes full responsibility for the content of the published article.

\end{document}